\documentclass[a4paper,12pt]{amsart}
\usepackage{amsfonts}
\usepackage{amssymb}
\usepackage{ifthen}
\usepackage{graphicx}
\nonstopmode \numberwithin{equation}{section}
\usepackage[usenames]{color}
\newtheorem{thm}{Theorem}
\newtheorem{lem}{Lemma}
\newtheorem{cor}{Corollary}

\newtheorem{cl}{Claim}
\newtheorem{ca}{Case}
\newtheorem{sca}{Subcase}
\newtheorem{scl}{Subclaim}
\newtheorem{conj}[equation]{Conjecture}

\theoremstyle{definition}
\newtheorem{defn}{Definition}

\newtheorem{op}[equation]{Open Problem}
\newtheorem{ques}[equation]{Question}
\newtheorem{rem}{Remark}[section]
\newtheorem{exam}[equation]{Example}

\newcounter {own}
\def\theown {\thesection       .\arabic{own}}

\newenvironment{pf}[1][]{%
 \vskip 3mm
 \noindent
 \ifthenelse{\equal{#1}{}}%
  {{\slshape Proof. }}%
  {{\slshape #1.} }%
 }%
{\qed\bigskip}

\newcounter{alphabet}
\newcounter{tmp}
\newenvironment{Thm}[1][]{\refstepcounter{alphabet}%
\bigskip%
\noindent%
{\bf Theorem \Alph{alphabet}}%
\ifthenelse{\equal{#1}{}}{}{ (#1)}%
{\bf .} \itshape}{\vskip 8pt}

\makeatletter
\newcommand{\Ref}[1]{\@ifundefined{r@#1}{}{\setcounter{tmp}{\ref{#1}}\Alph{tmp}}}
\makeatother

\newcommand{\IR}{{\mathbb R}}

\newcommand{\diam}{{\operatorname{diam}}}

\def\be{\begin{equation}}
\def\ee{\end{equation}}

\newcommand{\bee}{\begin{enumerate}}
\newcommand{\eee}{\end{enumerate}}

\newcommand{\blem}{\begin{lem}}
\newcommand{\elem}{\end{lem}}
\newcommand{\bthm}{\begin{thm}}
\newcommand{\ethm}{\end{thm}}
\newcommand{\bcor}{\begin{cor}}
\newcommand{\ecor}{\end{cor}}
\newcommand{\beg}{\begin{exam}}
\newcommand{\eeg}{\end{exam}}
\newcommand{\begs}{\begin{examples}}
\newcommand{\eegs}{\end{examples}}
\newcommand{\bdefe}{\begin{defn}}
\newcommand{\edefe}{\end{defn}}
\newcommand{\bprob}{\begin{prob}}
\newcommand{\eprob}{\end{prob}}
\newcommand{\bques}{\begin{ques}}
\newcommand{\eques}{\end{ques}}
\newcommand{\bei}{\begin{itemize}}
\newcommand{\eei}{\end{itemize}}
\newcommand{\bcon}{\begin{conj}}
\newcommand{\econ}{\end{conj}}
\newcommand{\bop}{\begin{op}}
\newcommand{\eop}{\end{op}}

\newcommand{\bca}{\begin{ca}}
\newcommand{\eca}{\end{ca}}
\newcommand{\bsca}{\begin{sca}}
\newcommand{\esca}{\end{sca}}

\newcommand{\bcl}{\begin{cl}}
\newcommand{\ecl}{\end{cl}}

\newcommand{\bscl}{\begin{scl}}
\newcommand{\escl}{\end{scl}}

\newcommand{\bcons}{\begin{conjs}}
\newcommand{\econs}{\end{conjs}}
\newcommand{\bprop}{\begin{propo}}
\newcommand{\eprop}{\end{propo}}
\newcommand{\br}{\begin{rem}}
\newcommand{\er}{\end{rem}}
\newcommand{\brs}{\begin{rems}}
\newcommand{\ers}{\end{rems}}
\newcommand{\bo}{\begin{obser}}
\newcommand{\eo}{\end{obser}}
\newcommand{\bos}{\begin{obsers}}
\newcommand{\eos}{\end{obsers}}
\newcommand{\bpf}{\begin{pf}}
\newcommand{\epf}{\end{pf}}
\newcommand{\ba}{\begin{array}}
\newcommand{\ea}{\end{array}}
\newcommand{\beq}{\begin{eqnarray}}
\newcommand{\beqq}{\begin{eqnarray*}}
\newcommand{\eeq}{\end{eqnarray}}
\newcommand{\eeqq}{\end{eqnarray*}}

\newcounter{minutes}
\divide\time by 60
\newcounter{hours}
\multiply\time by 60 \addtocounter{minutes}{-\time}

\begin{document}

\bibliographystyle{amsplain}
\title{Inner uniform domains and the Apollonian inner metric}
\thanks{$^\dagger$ File:~\jobname .tex,
          printed: \number\year-\number\month-\number\day,
          \thehours.\ifnum\theminutes<10{0}\fi\theminutes}
\author{Yaxiang Li}
\address{Yaxiang Li, College of Science,
Central South University of
Forestry and Technology, Changsha,  Hunan 410004, People's Republic
of China} \email{yaxiangli@163.com}
%
%
\author{Xiantao Wang
}
\address{Xiantao Wang, Department of Mathematics,
Hunan Normal University, Changsha,  Hunan 410081, People's Republic
of China} \email{xtwang@hunnu.edu.cn}

\date{}

\begin{abstract}
In this paper, we characterize inner uniform domains in $\IR^n$ in
terms of Apollonian inner metric and the metric $j'_D$ when $D$ are
Apollonian. As an application, a new characterization for $A$-uniform
domains is obtained.
\end{abstract}

\thanks{The research was partly supported by NSF of
China (No. 11071063) and Hunan Provincial Innovation Foundation For
Postgraduate.}

\maketitle\pagestyle{myheadings} \markboth{}{Inner uniform domains
and the Apollonian inner metric}

\section{Introduction and main results}\label{sec-1}

Throughout the paper, we assume that $D$ is a proper subdomain of
the Euclidean $n$-space $\mathbb{R}^n$, $n\geq 2$,
$[x,y]$ denotes the closed segment between $x$ and $y$, and
$B^n(x,r)$ stands for the open ball centered at $x$ with radius $r>0$,
i.e., $B^n(x,r)=\{y\in\mathbb{R}^n: |y-x|< r\}$. In particular,  we
use $\mathbb{B}^n$ to denote the unit ball $B^n(0,1)$. For $x, y\in
D$, the Apollonian distance is defined by
$$\alpha_D(x,y)= \sup_{a,b\in\partial D}\Big\{\log \frac{|a-x||b-y|}{|a-y||b-x|}\Big\},$$
where $\partial D$ means the boundary of $D$.
If one of $a$, $b$ equals to $\infty$, we understand that
$\frac{|\infty-x|}{|\infty-y|}=1.$ We note that this metric is
invariant under M\"{o}bius transformations and equals the hyperbolic
distance in balls and half spaces (cf. \cite{Bea}). It is in fact a
metric if and only if the complement of $D$ is not contained in a
hyperplane as was noted in \cite[Theorem 1.1]{Bea} (see also
\cite{H1}). In this paper, these domains are called to be {\it
Apollonian}.
 This metric was introduced in
\cite{Bea} and considered in \cite{Bar,GH,H,H1,H2,H3,
H4, Rh, Se}.

Let $\gamma: [0,1]\rightarrow D$ be a path, i.e., a continuous
function. If $d$ is a metric in $D$, then the $d$-length of $\gamma$
is defined by $$d(\gamma)=\sup \Big\{
\sum_{i=0}^{k-1}d(\gamma(t_i),\gamma(t_{i+1}))\Big\},$$ where the
supremum is taken over all $k<\infty$ and all sequences $\{t_i\}$
satisfying $0=t_0<t_1<\dots<t_k=1$. All the paths in this paper are
assumed to be rectifiable, that is, they have the finite Euclidean
arc length. The inner metric of $d$ is defined by the formula
$$\widetilde{d}(x,y)= \inf_{\gamma} \{d(\gamma)\},$$ where the
infimum is taken over all paths connecting $x$ and $y$ in $D$.
Particularly, we use $\widetilde{\alpha}_D$ to denote the inner
metric of the Apollonian metric $\alpha_D$ and call it the {\it
Apollonian inner metric}. Also we use $\lambda_D(x, y)$ to denote
$\widetilde{d}(x,y)$ when $d(\gamma)$ is the Euclidean arc
length.

In \cite[theorem 1.2]{H}, H\"{a}st\"{o} proved that
$\widetilde{\alpha}_D$ is a metric if and only if the complement of
$D$ is not contained in an $(n-2)$-dimensional hyperplane in
$\mathbb{R}^n$.
Further, in
\cite{H}, H\"{a}st\"{o} showed

\begin{Thm}\label{thmB}$($\cite[Theorem 1.5]{H}$)$ Let $D$ be
Apollonian. Then for $x$, $y\in D$, there exists a path $\gamma$ in
$D$ connecting $x$ and $y$ such that
$$\alpha_D(\gamma)=\widetilde{\alpha}_D(x,y).$$
\end{Thm}
And further, in \cite{HPWS}, the authors got the following.

\begin{Thm}\label{Lemc1}{\rm (\cite[Lemma 2.4]{HPWS})} Let $x,y\in
D$ and let $\gamma\subset D$ be a path such that
$\widetilde{\alpha}_D(x,y)=\alpha_D(\gamma).$ Then for each
$z$, $w\in \gamma$, we have
$$\widetilde{\alpha}_D(z,w)=\alpha_D(\gamma[z,w]),$$
where
$\gamma[z,w]$ denotes the part of $\gamma$ between $z$ and
$w$.\end{Thm}



\bdefe\label{def1} A domain $D$ is called {\it inner $c$-uniform} provided there exists a positive constant $c$ such that each
pair of points $z_{1},z_{2}$ in $D$ can be joined by a rectifiable
arc $\gamma$ in $D$ satisfying $($cf. \cite{Vai1}$)$\begin{enumerate}

\item $\min\{\ell (\gamma [z_1, z])\; \ell (\gamma [z_2, z])\}\leq c\, d_D(z)$ for all $z\in \gamma$; and

\item $\ell(\gamma)\leq c\,\lambda_D(z_1,z_2)$.\end{enumerate}where
$d_D(z)$ denotes the distance from $z$ to the boundary $\partial D$
of $D$.

If $\lambda_D(z_1,z_2)$ is replaced by $|z_1-z_2|$ in Definition
\ref{def1}, then $D$ is said to be $c$-{\it uniform}.\edefe
Obviously,
uniformity implies inner uniformity.

\bdefe\label{defe}A domain $D$ is called to be a {\it $c$-John
domain} provided there exists a positive constant $c$ such that each
pair of points $z_{1},z_{2}$ in $D$ can be joined by a rectifiable
arc $\gamma$ in $D$ satisfying $($cf. \cite{NV}$)$
$$\min\{\diam(\gamma[z_1,w]),\diam(\gamma[w,z_2])\}\leq c d_D(w).$$
\edefe

In \cite{Vai1}, V\"ais\"al\"a showed the following two theorems.

\begin{Thm}\label{ThmE} $($\cite[Theorem 3.3 and Theorem 3.4]{Vai1}$)$ Suppose that $D\subset R^n$ is
an inner $c$-uniform domain. Then for $x,y \in D$, we have
$$\lambda_D(x,y)\leq \nu_1\varrho_D(x,y),$$  where $\nu_1\geq 6c$ is a constant depending on $c$ and $n$, and $\varrho_D(x,y)$ denotes the inner diameter metric,
defined by $$\varrho_D(x,y)=\inf_{\gamma} \{\diam(\gamma)\}$$ over all arcs
$\gamma$ joining $x$ and $y$ in $D$.
\end{Thm}

\begin{Thm}\label{lem-vai1}$($\cite[Theorem 3.11]{Vai1}$)$ For a
domain $D\subset \mathbb{R}^n$, the following conditions are
quantitatively equivalent:\begin{enumerate}

\item $D$ is inner $c$-uniform.
\item Each pair of points $z_1,z_2\in D$ can be joined by an arc
$\gamma$ such that for $w\in \gamma$,
$$\min\{\diam(\gamma[z_1,w]),\; \diam(\gamma[z_2,w])\}\leq \nu_2 d_D(w)\;\;\mbox{and}
\;\;\diam(\gamma)\leq \nu_2\varrho_D(z_1,z_2),$$\end{enumerate}
where the constants $c$ and $\nu_2$ depend on each other and $n$.
\end{Thm}

Let $D$ be a domain and $x,y\in D$. We write
$$j_D(x,y)=\log\Big(1+\frac{|x-y|}{\min\{d_D(x),d_D(y)\}}\Big).$$ Kim
\cite{K06} $($see also \cite{Vai1}$)$ introduced the following
version of the $j$-metric:
$$j'_{D}(x,y)=\log\Big(1+\frac{\varrho_D(x,y)}{\min\{d_D(x),d_D(y)\}}\Big),$$
 and the quasihyperbolic metric \cite{GP} is defined by
$$k_D(x,y)=\inf_{\gamma}\int_{\gamma}\frac{|dz|}{d_D(z)},$$
where the infimum is taken over all paths $\gamma$ joining $x$ and
$y$ in $D$.

We easily know from the proof of
\cite[Lemma 2.2]{Vai6-0} that for $x$, $y\in D$,  \beq\label{equa}j_D(x,y)\leq
j'_D(x,y)\leq k_D(x,y).\eeq
Further, we have

\begin{Thm}\label{lemA} For $x,y\in D$, the following hold true.
\begin{enumerate}


\item$($\cite[Corollary 3.2]{Bea}$)$ \quad $\big|\log\frac{d_D(x)}{d_D(y)}\big|\leq\alpha_D(x,y)\leq 2j_D(x,y)$;

\item$($\cite[Lemma 5.3]{H1}$)$\quad$\widetilde{j}_D(x,y)=k_D(x,y);$
\item $($\cite[Corollary 5.4]{H1}$)$ \quad $\widetilde{\alpha}_D(x,y)\leq 2
  k_D(x,y)$.\end{enumerate}\end{Thm}

In \cite{Geo}, Gehring and Osgood got a characterization of uniform
domains in terms of $k_D$ and $j_D$.

\begin{Thm}\label{lemA1}$($\cite[Corollary 1]{Geo}$)$ A domain $D$ is
$\mu$-uniform if and only if there exists a constant $\mu_1$ such
that $$k_D(z_1,z_2)\leq \mu_1 j_D(z_1,z_2)$$ for all $z_1,z_2\in D$,
where the constants $\mu$ and $\mu_1$ depend only on each
other.\end{Thm}

As a matter of fact, the above inequality appeared in \cite{Geo} in
a form with an additive constant on the right hand side: it was
shown by Vuorinen \cite[2.50]{Vu2} that the additive constant can be
chosen to be $0$. Moreover, in \cite{HPWS}, the authors proved the
following.

\begin{Thm}\label{thmA}$($\cite[Theorem 1.2]{HPWS}$)$ A domain $D\subset
\mathbb{R}^n$ is $\mu$-uniform if and only if there exists a
constant $\mu_2$ such that $\widetilde{\alpha}_D(x,y)\leq \mu_2
j_D(x,y)$ for any $x,y\in D,$ where the constants $\mu$ and $\mu_2$
depend only on each other.
\end{Thm}

See \cite{BHK,FW,Geo,HPWS,K,Martio-80,MS, Vai1, Vai,Vai6} for more
details on uniform domains and inner uniform domains.

 By Theorem \Ref{thmA}, one may ask that if we can
characterize inner uniform domains in terms of
$\widetilde{\alpha}_D$ and $j'_D$.
 The main aim of this paper is to consider this problem. Our result
 shows that the answer to this problem is affirmative. Combining with \cite[Theorem 2.1]{K} and Theorem \Ref{lem-vai1}, we
state our result in the following form.

\begin{thm}\label{thm1} Let $D$ be a proper subdomain of
$\mathbb{R}^n.$ If $D$ is Apollonian, then the followings are quantitatively equivalent.\begin{enumerate}
\item\label{thm1-1} $D$ is an inner $c$-uniform domain;

\item\label{thm1-2} There exists a constant $c_1$ such that $$k_D(x,y)\leq c_1 j'_D(x,y) \;\;\;\;\forall x,y\in D;$$
\item\label{thm1-3} There exists a constant $c_2$ such that \beq\label{thm-eq}\widetilde{\alpha}_D(x,y)\leq c_2 j'_D(x,y) \;\;\;\;\forall x,y\in
D;\eeq
\item\label{thm1-4} Each pair of points $x,y\in D$ can be joined by
an arc $\gamma$ such that for $w\in\gamma,$
$\min\{\diam(\gamma[x,w]),\diam(\gamma[y,w])\}\leq c_3 d_D(w)$ and
$\diam(\gamma)\leq c_3 \varrho_D(x,y)$,
\end{enumerate} where $c$, $c_1$, $c_2$ and $c_3$ are constants greater than $1$, and depend on each other and $n$.
\end{thm}

In \cite{H1}, H\"ast\"o proved the following result.

\begin{Thm}\label{Hasto}$($\cite[Proposition 6.6]{H1}$)$ Let $D\subset
\mathbb{R}^n$ be a domain. The following conditions are quantitatively
equivalent:\begin{enumerate}
\item $D$ is $A$-uniform with coefficient $K$, that is, there exist some constant
$K$ such that for $x, y\in D$,  $k_D(x, y)\leq K \alpha_D(x, y);$

\item $D$ is $\mu$-uniform and has the comparison property with some constant $L$;

\item $D$ is $\mu_3$-quasi-isotropic and $\widetilde{\alpha}_D\leq \mu_4 \alpha_D$,\end{enumerate}
where the constants $K$, $L$, $\mu$, $\mu_3$ and $\mu_4$ depend only
on each other.\end{Thm} Here we say that a domain $D\subset
\mathbb{R}^n$ has the {\it comparison property} if there exists a
constant $L$ such that $$j_D/L\leq \alpha_D \leq 2j_D,$$ and $D$ is
$\mu_3$-quasi-isotropic if
$$\limsup_{r\rightarrow 0}\frac{\sup\{\alpha_D(x,z): |x-z|=r\}}{\inf\{\alpha_D(x,y): |x-y|=r\}}\leq \mu_3$$
for every $x\in D$ $($See \cite{H1}$)$.

As an application of Theorem \ref{thm1}, we get a new characterization for
 $A$-uniform domains.

\begin{cor}\label{cor} Let $D\subset
\mathbb{R}^n$ be an Apollonian domain. The following conditions are
quantitatively equivalent:\begin{enumerate}
\item\label{cor1} $D$ is $A$-uniform with coefficient $K$;
\item\label{cor2} $D$ is $c$-inner uniform and
$j'_D(x,y)\leq \mu_5 \alpha_D(x,y)$ for all $x,y\in
D$,\end{enumerate} where the constants $c$, $K$ and $\mu_5$ depend
on each other and $n$.
\end{cor}
In the
next section, we will prove Theorem \ref{thm1}  and
Corollary \ref{cor}.

\section{Proofs of Theorem \ref{thm1}  and
Corollary \ref{cor}}\label{sec-2}

\subsection{ Proof of Theorem \ref{thm1}}
The implication \eqref{thm1-1} $\Rightarrow$ \eqref{thm1-2} follows
from \cite[Theorem 2.1]{K} and Theorem \Ref{ThmE}, and Theorem
\Ref{lemA} shows that \eqref{thm1-2} $\Rightarrow$ \eqref{thm1-3} is
true. The implication \eqref{thm1-4}$\Rightarrow$ \eqref{thm1-1}
follows from Theorem \Ref{lem-vai1}.
Hence to finish the proof of Theorem \ref{thm1}, it remains only one implication
\eqref{thm1-3} $\Rightarrow$ \eqref{thm1-4} to be checked.

Suppose that the assertion
\eqref{thm1-3} in the theorem holds. To prove the truth of the assertion of \eqref{thm1-4} in the theorem,
 we let $x$, $y\in D$. Without loss of generality, we assume that $d_D(x)\leq
d_D(y)$.  We consider the case where $|x-y|< d_D(x)$ and the case
where  $|x-y|\geq d_D(x)$, separately.

\bca\label{case1} $|x-y|< d_D(x).$\eca

Let $\gamma=[x,y]$ be the Euclidean line segment joining $x$ and
$y$. Clearly, $\gamma\subset D$,

$$\diam (\gamma)=|x-y|=\varrho_D(x,y)$$ and
$$\min\{\diam(\gamma[x,w]),\;  \diam(\gamma[y,w])\}\leq d_D(w)$$ for $w\in\gamma$. Thus the assertion \eqref{thm1-4} in the theorem
is true in this case.

\bca\label{case2} $|x-y|\geq d_D(x).$\eca

By Theorem \Ref{thmB} there exists a path $\gamma\subset D$
connecting $x$ and $y$ such that
$$\widetilde{\alpha}_D(x,y)=\alpha_D(\gamma).$$ By
compactness we see that there is a point $z_0$ in $\gamma$ which is
the first point along the direction from $x$ to $y$ satisfying
$$d_D(z_0)=\sup_{w\in\gamma}\{d_D(w)\}.$$
Let $m\geq 0$ be the integer such that $$2^m d_D(x)\leq d_D(z_0)<
2^{m+1}d_D(x),$$ and let $x_0$ be the first point of $\gamma[x,z_0]$
from $x$ to $z_0$ with \beq\label{eq1-1'}d_D(x_0)=2^m d_D(x).\eeq
Then we have \beq\label{eq1-1}d_D(x_0)\leq d_D(z_0)< 2 d_D(x_0).\eeq
 Let $x_1=x$, and let $x_2, \dots, x_{m+1}$ be the points such that for each $i\in \{2, \dots, m+1\}$, $x_i$ denotes
 the first point in
 $\gamma[x,z_0]$ along the direction from $x$ to $z_0$ satisfying
 $$d_D(x_i)= 2^{i-1}d_D(x_1).$$ Apparently,
 $x_{m+1}=x_0.$ If $x_0\not= z_0$, we denote $z_0$ by $x_{m+2}$. By the choice of $x_i$, we know that for each
 $i\in\{1,2,\dots,m\}$, \beq\label{eq1-2'} d_D(x_{i+1})=2
 d_D(x_i),\eeq and so  \beq\label{eq1-2} \varrho_D(x_i,x_{i+1})\geq d_D(x_{i+1})-d_D(x_i)= d_D(x_i).\eeq
 For each $i\in\{1,2,\dots,m\}$ and $w\in \gamma[x_i,x_{i+1}]$, it easily follows that
\beq\label{eq1-3}d_D(w)\leq
 d_D(x_{i+1})= 2
 d_D(x_i).\eeq

Let $w_0$ be the first point of $\gamma$  along the direction from
$y$ to
 $x$ satisfying
 $$d_D(w_0)=\sup_{w\in\gamma}\{d_D(w)\}.$$
 Obviously, $d_D(w_0)=d_D(z_0)$. It is possible that $w_0=z_0$. A similar argument as above shows that there are points $\{y_j\}_{j=1}^{s+1}$ in
 $\gamma[y,w_0]$ such that for each $j\in\{1, \dots, s+1\}$, $y_j$ denotes the first point in $\gamma[y,w_0]$
 along the direction from $y$ to $w_0$ satisfying $$d_D(y_j)=2^{j-1}d_D(y_1),$$
 where $y_1=y$ and $d_D(y_{s+1})=2^sd_D(y_1)$. We also use $y_0$
 to denote $y_{s+1}$.
 If $y_0\not=w_0$, we use $y_{s+2}$ to
 denote $w_0.$

\begin{lem}\label{Lemma1} For each $i\in\{1,2,\dots,m\},$ we have
\begin{enumerate}
\item\label{cl1}$\diam(\gamma[x_i,x_{i+1}])\leq
b_1\varrho_D(x_i,x_{i+1})$ with $b_1= 24 c'_2$ and $c'_2=[c_2]+1.$
Here and in the following, $[\cdot]$ always denotes the greatest
integer part;
\item\label{cl2}$\varrho_D(x_i,x_{i+1})\leq b_2d_D(x_i)$ with $b_2= (1+b_1)^2;$
\item\label{cl3}$d_D(x_i)\leq b_3 d_D(w)$ for all $w\in \gamma$ with
$b_3=(1+b_2)^{\frac{c_2}{2}},$
\end{enumerate}
\noindent where $c_2$ is the same constant as in the inequality
\eqref{thm-eq}.\end{lem}
\bpf We now prove the first
assertion in the lemma.
 Suppose on the contrary that there is some $i\in\{1,\dots,m\}$
 satisfying
$$\diam(\gamma[x_i,x_{i+1}])>b_1 \varrho_D(x_i,x_{i+1}).$$ Let
$u_{i,1}=x_i$, and take the points
$u_{i,2},u_{i,3},\dots,u_{i,c'_2+1}$ in $\gamma$ such that for each
$t\in \{2, \dots, c'_2+1\}$, $u_{i,t}$ is the first point of
$\gamma$ from $x_i$ to $x_{i+1}$ satisfying
$$|x_i-u_{i,t}|=6(t-1)\varrho_D(x_i,x_{i+1}).$$ Then for each $t\in
\{1,\dots,c'_2\}$, we have \beq\label{eq1-4}|u_{i,t}-u_{i,t+1}|\geq
|x_i-u_{i,t+1}|-|u_{i,t}-x_i|\geq 6 \varrho_D(x_i,x_{i+1}).\eeq Let
$p\in \partial D$ be such that $d_D(u_{i,t+1})=|u_{i,t+1}-p|$. Then
\eqref{eq1-2}, \eqref{eq1-3} and \eqref{eq1-4} yield
\beq\label{eq1-5}|u_{i,t}-p|&\geq&
|u_{i,t}-u_{i,t+1}|-d_D(u_{i,t+1})\\\nonumber&\geq&
 6\varrho_D(x_i,x_{i+1})-2 d_D(x_i)\\\nonumber&\geq& 2\varrho_D(x_i,x_{i+1})+
2d_D(x_{i}).\eeq Similarly, for $q\in \partial D$ with
$d_D(u_{i,t})=|u_{i,t}-q|$, we get
\beq\label{eq1-6}|u_{i,t+1}-q|\geq 2\varrho_D(x_i,x_{i+1})
+2d_D(x_i).\eeq Thus we infer from \eqref{eq1-3}, \eqref{eq1-5} and
\eqref{eq1-6} that
$$\alpha_D(u_{i,t},u_{i,t+1})\geq \log\left(\frac{|u_{i,t}-p|}{d_D(u_{i,t+1})}\frac{|u_{i,t+1}-q|}{d_D(u_{i,t})}\right)\geq
2\log\left(1+\frac{\varrho_D(x_i,x_{i+1})}{d_D(x_i)}\right),$$ which
together with Theorem \Ref{Lemc1} show that
\begin{eqnarray*} \widetilde{\alpha}_D(x_i,x_{i+1})&=&
\alpha_D(\gamma[x_i,x_{i+1}])\\&\geq&
\sum_{t=1}^{c'_2}\alpha_D(u_{i,t},u_{i,t+1})\\&\geq&2c'_2\log\left(1+\frac{\varrho_D(x_i,x_{i+1})}{d_D(x_i)}\right)\\
&\geq&
2c_2\log\left(1+\frac{\varrho_D(x_i,x_{i+1})}{d_D(x_i)}\right)\\&=&2c_2
j'_D(x_i,x_{i+1}),\end{eqnarray*} which contradicts with
\eqref{thm-eq}. Hence \eqref{cl1} is true.

Then we come to prove the second assertion.  Suppose on the contrary
that there is some $i\in\{1,2,\dots,m\}$
 satisfying \beq\label{eq11}\varrho_D(x_i,x_{i+1})>
b_2d_D(x_i).\eeq

Obviously, there exists some point $v\in \gamma[x_i,x_{i+1}]$ such
that $|x_i-v|\geq \frac{1}{2}\varrho_D(x_i,x_{i+1}).$ We let
$v_{i,1}=x_i$, and let $v_{i,2}, \dots, v_{i,\frac{b_1}{12}+1}$ be
the points in $\gamma$ such that for each $h\in \{2, \cdots,
\frac{b_1}{12}+1\}$, $v_{i,h}$  is the first point of $\gamma$ from
$x_i$ to $x_{i+1}$ satisfying
$$|x_i-v_{i,h}|=\frac{6(h-1)}{b_1}\varrho_D(x_i,x_{i+1}).$$ Then
\beq\label{eq1-7}|v_{i,h}-v_{i,h+1}|\geq
|v_{i,h+1}-x_i|-|v_{i,h}-x_i|\geq
\frac{6}{b_1}\varrho_D(x_i,x_{i+1}).\eeq Let $p\in
\partial D$ satisfy $d_D(v_{i,h+1})=|v_{i,h+1}-p|$. Then it follows from
\eqref{eq1-2}, \eqref{eq1-3}, \eqref{eq11} and \eqref{eq1-7} that
\beq\label{eq1-8}|v_{i,h}-p|&\geq&
|v_{i,h}-v_{i,h+1}|-d_D(v_{i,h+1})\\\nonumber&\geq&
\frac{6}{b_1}\varrho_D(x_i,x_{i+1})-2 d_D(x_i)
\\\nonumber&>& \frac{2}{b_1}\varrho_D(x_i,x_{i+1})+
2d_D(x_{i}).\eeq Similarly, for $q\in \partial D$ with
$d_D(v_{i,h})=|v_{i,h}-q|$, we know
\beq\label{eq1-9}|v_{i,h+1}-q|\geq
\frac{2}{b_1}\varrho_D(x_i,x_{i+1}) +2 d_D(x_i).\eeq Thus we infer
from \eqref{eq1-3}, \eqref{eq11}, \eqref{eq1-8} and \eqref{eq1-9}
that
\begin{eqnarray*}\alpha_D(v_{i,h},v_{i,h+1})&\geq&
 \log\left(\frac{|v_{i,h}-p|}{d_D(v_{i,h+1})}\frac{|v_{i,h+1}-q|}{d_D(v_{i,h})}\right)\\&\geq&
2\log\left(1+\frac{\varrho_D(x_i,x_{i+1})}{b_1d_D(x_i)}\right)\\&>&
\frac{12c_2}{b_1}\log\left(1+\frac{\varrho_D(x_i,x_{i+1})}{d_D(x_i)}\right).\end{eqnarray*}

Whence Theorem \Ref{Lemc1} yields
\begin{eqnarray*} \widetilde{\alpha}_D(x_i,x_{i+1})&=&
\alpha_D(\gamma[x_i,x_{i+1}])\\&\geq&
\sum_{h=1}^{\frac{b_1}{12}}\alpha_D(v_{i,h},v_{i,h+1})
\\
&>&
c_2\log\left(1+\frac{\varrho_D(x_i,x_{i+1})}{d_D(x_i)}\right)\\&=&c_2
j'_D(x_i,x_{i+1}),\end{eqnarray*} which is the desired
contradiction.

To finish the proof of Lemma \ref{Lemma1}, it remains to check
\eqref{cl3}. Let $w\in \gamma$. Then \eqref{cl2} in the lemma,
\eqref{thm-eq}, Theorems \Ref{Lemc1} and \Ref{lemA} lead to
\begin{eqnarray*}2\log
\frac{d_D(x_i)}{d_D(w)}&<&\alpha_D(x_i,w)+\alpha_D(w,x_{i+1})\\&\leq&\alpha_D(\gamma[x_i,x_{i+1}])
\\&=&\widetilde{\alpha}_D(x_i,x_{i+1})\\
&\leq&c_2\log\left(1+\frac{\varrho_D(x_i,x_{i+1})}{d_D(x_i)}\right)\\&\leq&
c_2\log(1+b_2),\end{eqnarray*} which shows that $$d_D(x_i)\leq
(1+b_2)^{\frac{c_2}{2}}d_D(w),$$ which shows that \eqref{cl3} is
true by taking $b_3=(1+b_2)^{\frac{c_2}{2}}$. Hence the proof of
Lemma \ref{Lemma1} is complete. \epf

Similarly, we know that \begin{lem}\label{Lemma2} For each
$j\in\{1,\dots,s\},$ we have
\begin{enumerate}
\item\label{cl11}$\diam(\gamma[y_j,y_{j+1}])\leq
b_1\varrho_D(y_j,y_{j+1})$;
\item\label{cl12}$\varrho_D(y_j,y_{j+1})\leq b_2d_D(y_j);$
\item\label{cl13}$d_D(y_j)\leq b_3 d_D(w)$ for all $w\in
\gamma$.\end{enumerate}\end{lem} \medskip

Suppose $x_0\neq y_0$. Then

\begin{lem}\label{Lemma3} For $w\in \gamma[x_0,y_0]$, we have
$$d_D(x_0)\leq b^2_3d_D(w)\;\; \mbox{and}\;\; \diam(\gamma[x_0,y_0])\leq b_4
d_D(w),$$ where $b_4=b_1b_2b^2_3.$\end{lem} \bpf We note by
\eqref{eq1-1} that \beq\label{eq1-10}\frac{1}{2}d_D(y_0)< d_D(x_0)<
2 d_D(y_0), \eeq and for $w\in \gamma[x_0,y_0]$, we have
\beq\label{eq1-11}d_D(w)<d_D(z_0)< 2 d_D(x_0).\eeq We prove this
lemma by considering the case where $\varrho_D(x_0,y_0)\geq
d_D(x_0)$ and the case where $\varrho_D(x_0,y_0)< d_D(x_0)$,
separately.


Suppose first that $\varrho_D(x_0,y_0)\geq d_D(x_0).$ Then by
\eqref{eq1-10}, \eqref{eq1-11} and a similar argument as in the
proof of Lemma \ref{Lemma1}, we get for each $w\in \gamma[x_0,y_0]$,
\beq\label{lw-1'}d_D(x_0)\leq b^2_3d_D(w)\eeq and \beq\label{lw-1}
\diam(\gamma[x_0,y_0])\leq b_4 d_D(w),\eeq where $b_4=b_1b_2b^2_3.$

Suppose next that $\varrho_D(x_0,y_0)< d_D(x_0).$ In this case, we
need the following claim. \bcl\label{wl-0} For $w\in
\gamma[x_0,y_0]$, we have $|w-x_0|\leq (3^{c_2}+1)d_D(x_0)$. \ecl
Obviously, to prove this claim, it suffices to consider the case
$|x_0-w|\geq 2 d_D(x_0)$. Let $z\in \partial D$ satisfy
$|x_0-z|=d_D(x_0).$ Then it follows from \eqref{thm-eq},
\eqref{eq1-10}, Theorems \Ref{Lemc1} and \Ref{lemA} that
\beq\label{equ00}
\log\left(\frac{|w-x_0|}{d_D(x_0)}-1\right)&\leq&\log
\frac{|z-w|}{d_D(x_0)} \\\nonumber&\leq& \alpha_D(x_0,
w)\\\nonumber&\leq&
\alpha_D(\gamma[x_0,y_0])\\\nonumber&=&\widetilde{\alpha}_D(x_0,y_0)\\\nonumber&\leq&
c_2
\log\Big(1+\frac{\varrho_D(x_0,y_0)}{\min\{d_D(x_0),d_D(y_0)\}}\Big)\\\nonumber&\leq&
c_2 \log 3, \eeq from which the  claim easily follows. \medskip

By Claim \ref{wl-0}, we get
\beq\label{eq1}\diam(\gamma[x_0,y_0]) \leq
2(3^{c_2}+1) d_D(x_0).\eeq

Moreover, by Theorem \Ref{lemA} and a similar argument as in
\eqref{equ00}, we also have
$$\log \frac{d_D(x_0)}{d_D(w)}\leq \alpha_D(x_0, w)\leq c_2 \log
3,$$ and so \beq\label{lw-9} d_D(x_0)\leq 3^{c_2}d_D(w),\eeq which
together with \eqref{eq1} show that
\beq\label{lw-2}\diam(\gamma[x_0,y_0])\leq 2(3^{c_2}+1)^2d_D(w)\leq
b_4 d_D(w).\eeq The inequalities \eqref{lw-1'}, \eqref{lw-1},
\eqref{lw-9} and \eqref{lw-2} imply that the lemma is true.\epf

Now we come to prove that the first part of \eqref{thm1-4} in
Theorem \ref{thm1} holds with constant $2b_4$, i.e., for $w\in
\gamma$, \beq\label{lw-4}
\min\{\diam(\gamma[x,w]),\diam(\gamma[y,w])\}\leq 2b_4 d_D(w).\eeq

Let $w\in \gamma$. We divide the discussions into three cases.

\bca\label{ca-3} $w\in\gamma[x,x_0].$\eca Clearly, there exists an
integer $k\in\{1,2,\dots,m\}$ such that $w\in \gamma[x_k,x_{k+1}].$
By Lemma \ref{Lemma1} and \eqref{eq1-2'} we have
\beq\label{lw-8}\diam(\gamma[x,w])&\leq&
\sum_{i=1}^{k}\diam(\gamma[x_i,x_{i+1}])\leq
b_1b_2\sum_{i=1}^{k}d_D(x_i)\\ \nonumber &\leq& 2b_1b_2d_D(x_k)\leq
2b_1b_2b_3d_D(w)< b_4d_D(w). \eeq

\bca\label{ca-4} $w\in\gamma[y,y_0].$\eca By Lemma \ref{Lemma2}, we
see from a similar argument as in the proof of Case \ref{ca-3} that

\beq\label{lw-7}\diam(\gamma[y,w])< b_4d_D(w). \eeq \bca\label{ca-5}
If $w\in \gamma[x_0,y_0].$\eca It follows from Lemmas \ref{Lemma1}
and \ref{Lemma3} that \beq\label{lw-6}\diam(\gamma[x,w])&\leq&
\sum_{i=1}^{m}\diam(\gamma[x_i,x_{i+1}])+\diam(\gamma[x_0,y_0])\\
\nonumber &\leq& b_1b_2\sum_{i=1}^{m}d_D(x_i)+b_4d_D(w)\\ \nonumber
&\leq&b_1b_2d_D(x_{m+1})+b_4d_D(w)\\ \nonumber &\leq&2b_4d_D(w).
\eeq The proof for \eqref{lw-4} easily follows from the combination
of \eqref{lw-8}, \eqref{lw-7} and \eqref{lw-6}.

Next we prove the second part of \eqref{thm1-4} in Theorem
\ref{thm1} with constant $b_6=2^{b_5+2} b_4$, where $b_5=c_2\log_2
(8b_4+1)+1$, i.e., \beq\label{lw-5}\diam(\gamma[x,y])\leq
b_6\varrho_D(x,y).\eeq

We first prove a lemma.
 \begin{lem}\label{Lemma4}
$\varrho_D(x,y)\geq 2^{m-b_5}d_D(x)$.\end{lem} \bpf If $m\leq b_5$,
then it is obvious from the assumption ``$|x-y|\geq d_D(x)$". So we
assume that $m> b_5$. In this case, we prove the lemma by
contradiction. Suppose that \beq\label{eq-a}\varrho_D(x,y)<
2^{m-b_5}d_D(x).\eeq Then \beq\label{eq-b}d_D(y)&\leq&
\varrho_D(x,y)+d_D(x)\\\nonumber&<&
(2^{m-b_5}+1)d_D(x)\\\nonumber&\leq&
\frac{2^m}{(8b_4+1)^{c_2}}d_D(x)\\\nonumber&<&
\frac{2^m}{(8b_4+1)^{\frac{c_2}{2}}}d_D(x).\eeq

By \eqref{eq1-1'} and \eqref{eq1-10} we have \beq\label{lw-3}d_D(y_0) \geq
\frac{1}{2}d_D(x_0)=2^{m-1}d_D(x)>
\frac{2^m}{(8b_4+1)^{\frac{c_2}{2}}}d_D(x),\eeq then we obtain from
 \eqref{eq-b}, \eqref{lw-3} and the easy fact
``$\frac{2^m}{(8b_4+1)^{\frac{c_2}{2}}}=2^{m-\frac{c_2}{2}\log_2(8b_4+1)}>
1$" that
there exist $w_1\in\gamma[x, x_{0}]$ and $w_2\in \gamma[y,y_{0}]$
such that
 \beq\label{eq1-01}d_D(w_1)=d_D(w_2)=\frac{2^m}{(8b_4+1)^{\frac{c_2}{2}}}d_D(x).\eeq
On one hand, we obtain from \eqref{thm-eq}, \eqref{eq1-1'},
\eqref{eq1-01}, Theorems \Ref{Lemc1} and \Ref{lemA} that

\begin{eqnarray*}c_2\log\left(1+\frac{\varrho_D(w_1,w_2)}{d_D(w_1)}\right)&\geq&
\widetilde{\alpha}_D(w_1,w_2)\\&=&\alpha_D(\gamma[w_1,w_2])\\&\geq&
\alpha_D(w_1,x_{m+1})+\alpha_D(x_{m+1},w_2)\\&\geq&
2\log\frac{d_D(x_{m+1})}{d_D(w_1)}\\&\geq&
c_2\log(1+8b_4),\end{eqnarray*} which imply that
\beq\label{eq1-12'}\varrho_D(w_1,w_2)\geq 8b_4 d_D(w_1).\eeq On the
other hand, by \eqref{lw-8}, \eqref{lw-7}, \eqref{eq-a} and
\eqref{eq1-01} we obtain
\begin{eqnarray*}\varrho_D(w_1,w_2)&\leq&
\varrho_D(w_1,x)+\varrho_D(x,y)+\varrho_D(y,w_2)\\&\leq&
\diam(\gamma[x,w_1])+\varrho_D(x,y)+\diam(\gamma[y,w_2])\\&<& 2b_4
d_D(w_1) +\varrho_D(x,y)\\&\leq&(2b_4+
\frac{2^{m-b_5}(8b_4+1)^{\frac{c_2}{2}}}{2^m})d_D(w_1)\\&<&(2b_4+1)d_D(w_1),\end{eqnarray*}
which is contradict with \eqref{eq1-12'}. Hence the proof of the
lemma is complete.\epf
\medskip

Now we are ready to conclude the proof of \eqref{lw-5}. It follows
from \eqref{eq1-1'}, \eqref{eq1-10},  \eqref{lw-8}, \eqref{lw-7}, \eqref{lw-3},
Lemmas \ref{Lemma3} and \ref{Lemma4} that

\begin{eqnarray*}\diam(\gamma[x,y])&\leq&
\diam(\gamma[x,x_0])+\diam(\gamma[x_0,y_0])+\diam(\gamma[y_0,y])\\&\leq&
4 b_4 d_D(x_{m+1})\\&=& 2^{m+2} b_4 d_D(x)\\&\leq&
2^{b_5+2}b_4\varrho_D(x,y).
\end{eqnarray*}
Hence the proof of \eqref{thm1-4} of Theorem \ref{thm1} is complete
by taking $c_3=2^{b_5+2} b_4.$\qed

\subsection{ The proof of Corollary \ref{cor}}

\eqref{cor1}$\Rightarrow$ \eqref{cor2}. Suppose \eqref{cor1} holds.
Then by Theorem \Ref{Hasto}, we see that $D$ is a $\mu$-uniform, so
it is obvious inner uniform.

For the second part of \eqref{cor2}, we can obtain easily from
\eqref{equa} and the definition of $A$-uniform $($see Theorem
\Ref{Hasto}$)$. Hence \eqref{cor2} is true.

\eqref{cor2} $\Rightarrow$ \eqref{cor1}. Suppose \eqref{cor2} holds.
Then by Theorem \ref{thm1}, we know that for all $x,y\in D,$
$$k_D(x,y)\leq c_1j'_D(x,y) \leq c_1\mu_5 \alpha_D(x,y),$$ which shows
that $D$ is $A$-uniform with coefficient $K=c_1\mu_5.$ \qed

\bigskip
{\bf Acknowledgements}  The authors thank the referees who have made valuable comments on  this manuscript. 



\begin{thebibliography}{99}
\bibitem{Bar}  {\sc D. Barbilian}, Einordnung von Lobayschewskys Massenbestimmung in einer gewissen allgemeinen
Metrik der Jordansche Bereiche, \textit{Casopsis Mathematiky a
Fysiky}, {\bf 64} (1934-35), 182--183.

\bibitem{Bea}  {\sc A. F. Beardon}, The Apollonian metric of a domain in $\mathbb{R}^n$, in: Quasiconformal Mappings
and Analysis, edited by Peter Duren, Juha Heinonen, Brad Osgood, and
Bruce Palka (\textit{Springer-Verlag, New York}), 1998, 91--108.

\bibitem{BHK} {\sc M. Bonk, J. Heinonen and P. Koskela}, Uniformizing Gromov
hyperbolic domains, \textit{Asterisque} {\bf 270} (2001), 1--99.


\bibitem{FW}  {\sc F. W. Gehring}, Uniform domains and the ubiquitous
quasidisks, \textit{Jahresber. Deutsch. Math. Verein,} {\bf
89} (1987), 88--103.



\bibitem{GH}  {\sc F. W. Gehring and K. Hag}, The Apollonian metric and quasiconformal mappings,
\textit{In the tradition of Ahlfors and Bers $($Stony Brook,
Ny, 1998$)$ $($Irwin Kra and Bernard Maskit, eds$)$, Contemp. Math.
256, Amer. Math. Sor.}, Providence, RI, 2000, MR 2001c: 30043,
143--163.



\bibitem{Geo}  {\sc F. W. Gehring and B. G. Osgood}, Uniform domains
 and the
quasi-hyperbolic metric, \textit{J. Analyse Math.,} {\bf 36} (1979),
50--74.


\bibitem{GP}  {\sc F. W. Gehring and B. P. Palka}, Quasiconformally
homogeneous domains, \textit{J. Analyse Math.}, {\bf 30} (1976),
172--199.

\bibitem{H}  {\sc P. H\"{a}st\"{o}}, The Apollonian inner metric, \textit{Commun. Anal. Geom.}, {\bf 12} (2004),
927--947.


\bibitem{H1}  {\sc P. H\"{a}st\"{o}}, The Apollonian  metric: Uniformity and quasiconvexity, \textit{Ann. Acad. Sci. Fenn. Math.},
 {\bf 28} (2003),
385--414.
\bibitem{H2}  {\sc P. H\"{a}st\"{o}}, The Apollonian  metric: limits of the approximation and bilipschitz properties,
\textit{Abstr. Appl. Anal.},
 {\bf 2003} (2003),
1141--1158.


\bibitem{H3}  {\sc P. H\"{a}st\"{o}}, The Apollonian  metric: quasi-isotropy and Seittenranta's metris,
\textit{Comput. Methods Funct. Theory},
 {\bf 4} (2004),
249--273.

\bibitem{H4}  {\sc P. H\"{a}st\"{o}}, The Apollonian metric: the comparison property, bilipchitz mappings and thick sets,
\textit{J. Appl. Anal.},
 {\bf 12} (2006),
209--232.


\bibitem{HPWS}  {\sc M. Huang, S. Ponnusamy, X. Wang and S. K. Sahoo}, The Apollonian inner metric and uniform domain,
\textit{Math. Nachr.}, {\bf 283} (2010), 1277--1290.



\bibitem{K06}  {\sc K. W. Kim}, The quasihyperbolic metric and analogues of the Hardy-Littlewood property for $\alpha=0$
in uniformity John domains,  \textit{Bull. Korean. Math. Soc.}, {\bf
43} (2006), no 2, 395--410.

\bibitem{K}  {\sc K. W. Kim}, Inner uniform domains, the
quasihyperbolic metric and weak bloch functions,  \textit{Bull.
Korean. Math. Soc.}, {\bf 49} (2012), 11--24.





\bibitem{Martio-80}  {\sc O. Martio},
Definitions of uniform domains, \textit{Ann. Acad. Sci. Fenn. Ser. A
I Math.,} {\bf 5} (1980), 197--205.

\bibitem{MS}  {\sc O. Martio and J. Sarvas},
Injectivity theorems in plane and space, \textit{Ann. Acad. Sci.
Fenn. Ser. A I Math.,} {\bf 4} (1978), 383--401.

\bibitem{NV}  {\sc R. N\"{a}kki and J. V\"{a}is\"{a}l\"{a}}, John disks,
\textit{Expo. Math.} {\bf 9} (1991), 3--43.


\bibitem{Rh}  {\sc A. G. Rhodes}, A upper bound for the hyperbolic metric of convex domain,  \textit{Bull.
London. Math. Soc.}, {\bf 29} (1997), 592--594.

\bibitem{Se}  {\sc P. Seittenranta}, M\"obius-invariant metrics,  \textit{Math. Proc. Cambridge Philos. Soc.},
 {\bf 125} (1999), 511--533.


\bibitem{Vai1}  {\sc J. V\"{a}is\"{a}l\"{a}}, Relatively and inner
 uniform domains, \textit{Conformal Geom. Dyn.,} {\bf 2} (1998), 56--88.


\bibitem{Vai}  {\sc J. V\"{a}is\"{a}l\"{a}}, Uniform domains,
\textit{Tohoku Math. J.,} {\bf 40} (1988), 101--118.



\bibitem{Vai6-0}  {\sc J. V\"{a}is\"{a}l\"{a}}, Free quasiconformality
in Banach spaces I, \textit{Ann. Acad. Sci. Fenn. Ser. A I Math.,}
{\bf 15} (1990), 355-379.

\bibitem{Vai6}  {\sc J. V\"{a}is\"{a}l\"{a}}, Free quasiconformality
in Banach spaces II, \textit{Ann. Acad. Sci. Fenn. Ser. A I Math.,}
{\bf 16} (1991), 255-310.

\bibitem{Vu2}  {\sc M. Vuorinen,}
Conformal invariants and quasiregular mappings, {\em J. Anal.
Math.}~\textbf{45} (1985), 69--115.

\end{thebibliography}
\end{document}